\documentclass[12pt,seceqn,secthm,amsthm
]{elsart}
\journal{Appl. Comput. Harmon. Anal.}
\usepackage{amsmath,amssymb}
\newtheorem{theorem}{Theorem}[section]
\theoremstyle{plain}
\newtheorem{lemma}[theorem]{Lemma}
\theoremstyle{definition}
\newtheorem{definitions}[theorem]{Definitions}
\newtheorem*{acknowledgements}{Acknowledgements}
\newtheorem{remark}[theorem]{Remark}
\def\openone{\hbox{\upshape \small1\kern-3.3pt\normalsize1}}
\newcounter{saveenumi}
\newlength{\byperiodicity}
\begin{document}
\begin{frontmatter}

\title%
{Sample Paths in Wavelet Theory}%
\author
{Palle E. T. Jorgensen\corauthref{cor1}\thanksref{label1}}
\ead{jorgen@math.uiowa.edu}
\ead[url]{http://www.math.uiowa.edu/\symbol{126}jorgen/}
\address
{Department of Mathematics, The University of Iowa, 14 MacLean Hall, \mbox{Iowa City,}
IA 52242-1419, U.S.A.}
\corauth[cor1]{\textit{Tel:} (319) 335-0782, \textit{Fax:} (319) 335-0627}
\thanks[label1]{This material is based upon work supported by the U.S. National Science
Foundation under Grant No.\ DMS-0139473 (FRG)}%
\begin{abstract}
We consider a class of convergence questions for infinite products that arise in wavelet theory when the wavelet filters are more singular than is traditionally built into the assumptions. We establish pointwise convergence properties for the absolute square of the scaling functions. Our proofs are based on probabilistic tools.
\end{abstract}
\begin{keyword}
scaling identity
\sep
quadrature filter
\sep
low-pass
\sep
Kolmogorov
\sep
sampling
\sep
pointwise convergence
\renewcommand{\MSC}{{\par\leavevmode\hbox {\it 2000 MSC:\ }}}
\MSC
42C40
\sep
42A16
\sep
43A65
\sep
42A65
\end{keyword}%
\end{frontmatter}

\section{\label{IntWav}Introduction. Wavelets}

We will be primarily interested in the infinite products that
typically occur in wavelet analysis \cite[Ch.\ 5]{Dau92}; see also \cite{BeLe01}. As is
well known, the traditional approach to wavelet-scaling functions is based
on analysis of wavelet filters which are periodic functions in one or more
variables, and which have some degree of regularity (say, Lipschitz at low
frequencies). In addition, the filter functions must satisfy a variant of
the following two features: (i)~some quadrature condition, and (ii)~some
low-pass condition. But conditions (i)--(ii) are known to be incompatible
with several new classes of examples. For these wider classes of wavelets,
we show that a new \emph{uncertainty principle} for wavelet filters is at play.
Reflecting this uncertainty principle for wavelet filters, we note that
recent papers on wavelet sets, and on other frequency-localized wavelets
(see, e.g., \cite{BJMP04,BMM99,DaLa98,DLS97}), necessitate the consideration
of wavelet filters which are defined only a.e., and for which conditions
(i)--(ii) are more subtle. As a result, the standard deterministic methods
then do not immediately apply.

In this paper, we use a random walk model (see
\cite{DiFr99,Spi76}) 
in a new analysis of
wavelets,
but for a class of wavelets where the filters are more
singular than is the case in the standard references \cite{Dau92,Gro01,HeWe96,JMR01,Mal98,MeCo97,Wal02}. 
Our results apply also to iterative algorithms
for wavelet packets \cite{Wic93,Wic94}. In these problems solutions are
constructed from an algorithm%
\index{algorithm}
which relates a certain function $\varphi(t)$, $t\in\mathbb{R}$, called a
scaling function, to its scaled version $\varphi(Nt)$ where $N$ is a fixed
integer, $2$ or more, or an expansive integral $d$ by $d$ matrix $A$ if $t$ is
in $\mathbb{R}^{d}$.

For $d = 1$, the formula which determines $\varphi$ is 
\begin{equation}
\varphi\left(  t\right)  =N\sum_{k\in\mathbb{Z}}a_{k}\varphi\left(
Nt-k\right)  ,\qquad t\in\mathbb{R}. \label{eqInt.7}%
\end{equation}
It is called the \emph{scaling identity}.
Consider the filter function
$m\left(  x\right)  =\sum\limits_{\makebox[0pt]{\hss$\scriptstyle k\in\mathbb{Z}$\hss}}a_{k}e^{-i2\pi kx}$
for $z = \exp (-i 2 \pi x)$ in the circle $\mathbb{T}$, or the
$d$-torus $\mathbb{T}^d$. It is traditional to impose some degree of regularity, e.g.,
Lipschitz, or low-pass for $x = 0$; but we shall not do this. Instead we shall
work with a probabilistic notion which serves as a minimal restriction on
the filter $m\left(  x\right) $ guaranteeing the existence of
$L^2$ solutions to the scaling identity (\ref{eqInt.7}).
It is not at all clear \emph{a priori} that this identity should have
$L^{2}$ solutions, or even
solutions that are given by some kind of convergent algorithm. The
coefficients in (\ref{eqInt.7}) are called \emph{masking coefficients}, \emph{filters}, or
filter coefficients. For special values of the coefficients, it turns out that
there is a normalized $L^{2}$ solution $\varphi$, and that pointwise
convergence of a good approximation can be established in a meaningful way.
Moreover, we show that the convergence behavior is dictated by properties of
an associated transfer operator $R$, or Ruelle operator.

A main point in this paper is our suggestion of a different and
more versatile approach to how we impose low-pass conditions on the basic
wavelet filters. Our approach is probabilistic, and it is based directly on
the Ruelle operator $R = R_W$, and on a certain Perron--Frobenius eigenfunction
$h$ for $R$ (see Theorem \ref{CorIntSam.4}).  Our work is motivated in part by recent
considerations of frequency-localized wavelets, as pioneered in for example
in the papers \cite{BCMO95,BMM99,BaMe99,DaLa98,DLS97}. These papers suggest
an interesting context for wavelets that are necessarily frequency-%
localized. But at the same time, it is evident from this work that the
corresponding wavelet filters will then typically not satisfy the low-pass
properties that are otherwise known to hold for more traditional time- (or
space-) localized wavelets.

The expression on the right-hand side in equation (\ref{eqInt.7}) is also
called a \emph{subdivision} because of the function values $\varphi( Nt -k)$.
An iteration of the operations on the right hand side in (\ref{eqInt.7}) on
some initial function is called a \emph{cascade approximation}. This approximation is
one approach to the function $\varphi$, and the other is the infinite product
formula (\ref{eqInt.10bis}) for the Fourier transform of $\varphi$.

The relation between $\varphi$ and its scaled refinement is well understood
when we pass from the time domain to the frequency domain, via the Fourier
transform. In that case the relation is multiplicative, and involves a certain
periodic matrix function $m$, called the \emph{low-pass filter}. For further details,
see formulas (\ref{eqInt.8}) and (\ref{eqInt.10per}) below.

The study of the filters $m$ is part of signal processing (see, e.g.,
\cite{Jor03,StNg96}). But by a `mathematical miracle', they have become one of
the most useful tools in wavelet constructions; and at the same time, they
have pointed to a host of exciting applications of wavelet mathematics; see,
e.g., \cite{HJL04}. To get a path space measure for some of the wavelet
problems, we use the quadrature mirror properties (see
(\ref{eqInt.3}) below), or their generalizations,
which are assumed for the \emph{filter function} $m$, also called a frequency
response function. It is periodic, in one or several variables. In
$\mathbb{R}^{d}$, there is a variety of choices of a period lattice for the
problems at hand.

For a number of applications (see especially \cite{BMM99,BaMe99,DaLa98,DLS97}), we must consider vector versions of the scaling identity (\ref{eqInt.7}); and
in that case, the solution $\varphi$, the \emph{scaling function}, will be viewed as a
vector-valued function, i.e., a function from $\mathbb{R}^d$ into some Hilbert space,
typically finite-dimensional. In that case, the coefficients $a_k$ in (\ref{eqInt.7}) will
be matrix-valued; and the product on the right-hand side in (\ref{eqInt.7}) will be
matrix acting on vector. Following \cite{BMM99}, we then say that the initial
resolution subspace $V_0$ in $L^2(\mathbb{R}^d)$ has multiplicity. In the more
traditional multiresolution analysis (MRA) approach to wavelets, $\varphi$ is a
scalar function, and $V_0$ is the closed span of the $\mathbb{Z}$-translates of $\varphi$. The
case of multiplicity is called the generalized MRA (GMRA) approach. The
present considerations deal with the issue of passing from the filter
function (possibly matrix-valued) to the scaling function $\varphi$. However, the
formulation of our ideas in the scalar case may easily be modified to the
matrix/vector case; and for the sake of simplicity, our technical discussion
below will be presented in the scalar case. We leave to the reader the spelling
out of the generalization to GMRAs.

     The fact that there are solutions in $L^2(\mathbb{R}^d)$ is not at all obvious; see 
\cite{Dau92}. In application to images, the subspace $V_0$
(where $V_0:=\operatorname{cl\;span}\{\,\varphi(\,\cdot-k)\mid k\in\mathbb{Z}^d\,\}$)
may represent a certain 
resolution, and hence there is a choice involved, but we know by standard 
theory, see, e.g., \cite{Dau92}, that under apropriate conditions, such choices are 
possible. As a result there are extremely useful, and computationally 
efficient, wavelet bases in $L^2(\mathbb{R}^d)$. A resolution subspace $V_0$ within $L^2(\mathbb{R}^d)$
can be chosen to be arbitrarily fine: Finer resolutions correspond to larger 
subspaces.

    As noted for example in \cite{BrJo02b,DuJo05c}, a variant of the scaling 
equation is also used in computer graphics: 
there data is successively subdivided and the refined level of data is related 
to the previous level by prescribed masking coefficients. The latter 
coefficients in turn induce generating functions which are direct analogues of 
wavelet filters.

     One reason for the computational efficiency of wavelets lies in the fact 
that \emph{wavelet coefficients}%
\index{coefficient!wavelet}
in \emph{wavelet expansions} for functions in $V_0$ may be 
computed using matrix iteration, rather than by a direct computation of inner 
products: the latter would involve integration over $\mathbb{R}^d$, and hence be 
computationally inefficient, if feasible at all. The deeper reason 
for why we can compute wavelet coefficients using matrix iteration is an 
important connection to the subband filter%
\index{filter!subband}%
ing method from signal/image 
processing%
\index{signal processing}%
\index{image processing}
involving digital filters, down-sampling and up-sampling. In this 
setting filters may be realized as functions $m_0$ on a $d$-torus, e.g., 
quadrature mirror filters. 

    As emphasized for example in \cite{Jor05a} and \cite{Bri95}, because of down-sampling, the matrix 
iteration involved in the computation of wavelet coefficients involves so-%
called slanted Toeplitz matrices $F$ from signal processing. 
The slanted matrices $F$ are immediately available; they 
simply record the numbers (masking coefficient%
\index{coefficient!masking}%
s) from the $\varphi$-scaling 
equation.
These matrices further have
the computationally attractive property that the iterated powers $F^k$ become 
sucessively more sparse as $k$ increases, i.e., the matrix representation of 
$F^k$ has mostly zeros, and the non-zero terms have an especially attractive 
geometric configuration. In fact subband%
\index{subband}
signal processing yields a finite 
family, $F$, $G$, etc., of such slanted matrices; for example
(with $L$ for ``low frequency%
\index{frequency}%
\index{filter!low-pass}%
'' and $H$ for ``high frequency%
\index{frequency}%
\index{filter!high-pass}%
''):
\begin{gather*}
\varphi=\sum_{k}P_{k}\varphi(2\,\cdot\,-k),
\quad\psi=\sum_{k}Q_{k}\varphi(2\,\cdot\,-k),\label{Sphinewa}
\\
\left\{\begin{aligned}
F\colon &\quad y_{n}^{L}=\frac{1}{\sqrt{2}}\sum_{k}P_{k-2n}x_{k},\\
G\colon &\quad y_{n}^{H}=\frac{1}{\sqrt{2}}\sum_{k}Q_{k-2n}x_{k}.
\end{aligned}
\right.
\end{gather*}
The wavelet coefficients at 
scaling level $k$ of a numerical signal $s$ from $V_0$ are then simply the 
coordinates of $G F^k s$. By this we mean that a signal in $V_0$ is represented by 
a vector $s$ via a fixed choice of scaling function; see \cite{Dau92,BrJo02b,Jor03}.

Given some choice of scaling function $\varphi$, then under suitable conditions
(e.g., orthogonality or \emph{a priori} frame estimates) we may define the
operator $W\colon V_0\rightarrow\ell^2$
which links the resolution subspace
$V_0$ with the sequence space $\ell^2$
of signals $s=(s_k)$ as follows:
\[
W\left(\sum_{k\in\mathbb{Z}^d}s_k\,\varphi\left(\,\cdot-k\right)\right)=\left(s_k\right).
\]
Here $W$ becomes a well defined linear operator,
$W\colon V_0\rightarrow\ell^2$

Then 
the matrix product $G F^k$ is applied to $s$; and the matrices $G F^k$ get more 
slanted as $k$ increases.

    Our approach begins with the observation that the computational feature of 
this engineering device can be said to begin with an endomorphism $r_A$ of the 
$d$-torus%
\index{torus}
$\mathbb{T}^d = \mathbb{R}^d/\mathbb{Z}^d$\label{STnewa}, an endomorphism which results from simply passing 
matrix multiplication by $A$ on $\mathbb{R}^d$ to the quotient by $\mathbb{Z}^d$. It is then immediate 
that the inverse images 
$r_A^{-1}(x)$ are finite for all $x$ in $\mathbb{T}^d$, in fact $\# r_A^{-1}(x) = |\det A |$. From 
this we recover the scaling identity, and we note that the wavelet scaling 
equation is a special case of a more general identity known in computational 
fractal theory%
\index{fractal!theory@--- theory}
and in symbolic dynamics \cite{Jor04}. We show that wavelet algorithms and 
harmonic analysis naturally generalize to affine iterated function systems. 
Moreover, in this general context, we are able to build the ambient Hilbert 
space%
\index{space!Hilbert!fractal}%
s for a variety of dynamical systems which arise from the iterated 
dynamics of endomorphisms of compact spaces \cite{DuJo03}. 

      As a consequence, the fact that the ambient Hilbert space in the 
traditional wavelet setting is the more familiar $L^2(\mathbb{R}^d)$ 
is merely an artifact of the choice of filters $m_0$. As we further show, by 
enlarging the class of admissible filters, there is a variety of other ambient 
Hilbert spaces possible with corresponding wavelet expansions: the 
most notable are those which arise from iterated function systems (IFS) of 
fractal type, for example for the middle-third Cantor set, and scaling by $3$.

With examples, with theorems, and with graphics, we hope to bring these
threads to light in this little book: The journey from wavelets to fractals
via signal processing.

     More generally, there is a variety of other natural 
dynamical settings (affine IFSs) that invite the same computational approach. 
 
     The two most striking examples which admit such a harmonic analysis 
are perhaps complex dynamics and subshifts. Both will be worked out in 
detail. 
In the first case, consider a given rational function $r(z)$ of one complex 
variable. We then get an endomorphism $r$ acting on an associated Julia set 
$X$ in the complex%
\index{dynamics!complex}
plane $\mathbb{C}$ as follows: This endomorphism $r\colon X \to X$  
results by restriction to $X$ \cite{Bea91}. (Details: Recall that $X$ is by 
definition the complement of the points in $\mathbb{C}$ where the sequence of iterations 
$r^n$ is a normal family. Specifically, the Fatou set%
\index{Fatou!set@--- set}
$F$ of $r(z)$ is the 
largest open set in $\mathbb{C}$ where $r^n$ is a normal sequence of functions, and we 
let $X$ be the complement of $F$. Here $r^n$ denotes the $n$'th iteration of the 
rational function $r(z)$.) The induced endomorphism $r$ of $X$ is then simply 
the restriction to $X$ of $r(z)$.
If $r$ then denotes the 
resulting endomorphism, $r\colon X \to X$,
it is known \cite{DuJo04a} that $\# r^{-1}(x) = {}$degree of $r$, for every $x$ 
in $X$ (except for a finite set of singular points).

    In the second case, for a particular one-sided subshift, we may take $X$ as 
the corresponding state space, and again we have a naturally induced finite-to-%
one endomorphism of $X$ of geometric and computational significance.

    But in the general framework, there is not a natural candidate for the 
ambient Hilbert space. That is good in one sense, as it means that the subband 
filter%
\index{filter!subband}%
s $m_0$ which are feasible will constitute a richer family of functions on 
$X$.
 
      In all cases, the analysis is governed by a random-walk model with 
successive iterations where probabilities 
are assigned on the finite sets $\# r^{-1}(x)$ and are given by the function $W := | 
m_0 |^2 $. This leads to a \emph{transfer operator}%
\index{operator!transfer}
$R_W$\label{SRWnewa} which has 
features in common with the classical operator considered first by Perron and 
Frobenius for positive matrices, in particular it has a Perron--Frobenius 
eigenvalue, and positive Perron--Frobenius eigenvectors, one on the right, a 
function%
\index{eigenfunction!Perron Frobenius@Perron--Frobenius}%
, and one on the left, a measure; see \cite{Rue89}. 
This Perron--%
Frobenius measure, also sometimes called the Ruelle measure%
\index{measure!Ruelle}%
, is an essential 
ingredient for our construction of an ambient Hilbert space. All of this, we 
show, applies to a variety of examples, and as we show, has the more 
traditional wavelet setup as a special case, in fact the special case when 
the Ruelle measure on $\mathbb{T}^d$ is 
the Dirac mass%
\index{measure!Dirac}
corresponding to the point $0$ in $\mathbb{T}^d$ (additive notation) 
representing zero frequency in the signal processing setup.

    There are two more ingredients entering in our construction of the ambient 
Hilbert space: a path-space measure governed by the $W$-probabilities, and 
certain finite cycle%
\index{cycle}%
s for the endomorphism $r$.
For each $x$ in $X$, we consider 
paths by infinite iterated tracing back with $r^{-1}$ and recursively assigning 
probabilities with $W$. Hence we get a measure $P_x$\label{SPxnewa} on a space of paths for each 
$x$. These measures are in turn integrated in $x$ using the Ruelle measure on $X$. 
The resulting measure will now define the inner product in the ambient Hilbert 
space. 

Since the first question is to decide when $\varphi$ is in $L^{2}$, we iterate
and get an infinite matrix product involving the matrix $W:=m^{\ast}m$. Since
$W$ is positive semidefinite, we may create a positive path measure of a
random walk starting at $x$ in some period interval. In several dimensions,
$x$ starts in a fundamental domain $D$ for some fixed lattice, for example
$\mathbb{Z}^{d}$. The paths starting at $x$ arise by iteration of the inverse
branches of $x\rightarrow Ax\bmod\mathbb{Z}^{d}$. There are $N=|\det A|$
distinct branches. These $N$ branches may be viewed as endomorphisms of $D$.

In this paper, we will consider pointwise convergence of the
infinite product (\ref{eqInt.10bis}) below for the Fourier transform $\hat{\varphi}$ of the
scaling function $\varphi$. But the traditional low-pass/regularity considerations
for the filter $m(x)$ are more general. Once pointwise convergence (Theorem
\ref{CorIntSam.4}) is established, then the quadrature property for $m(x)$ will imply that
the scaling function $\varphi$ is automatically in $L^2(\mathbb{R})$.

We construct our random walks in a general framework which includes both
wavelets, wavelet packets, and some of the other more
classical problems. We further show how some of the classical questions may be
phrased and solved with the use of path space measures.

It is interesting to contrast our proposed approach with the more traditional
one used in wavelet analysis: see, e.g., \cite{Dau92}. Traditionally, some
kind of Lipschitz or Dini regularity condition must be assumed for the filter.
Then the corresponding infinite product may be made precise, and we can turn
to the question of when the wavelet generators are in $L^{2}(\mathbb{R}^{d})$;
see, e.g., \cite{Jor01a} and \cite{DuJo03}. As it turns out, both of these
issues have natural formulations, and solutions, in terms of the path space
measures. And the results allow a wider generality. In a variety of wavelet
questions for band-limited wavelets, the regularity conditions just aren't
satisfied for the filters $m$ that are dictated by the setting and the
applications; see, e.g., \cite{BJMP04}.

\section{\label{IntPat}Path space}

A well tested tool in analysis, and in mathematical physics, centers around
the application of path-space measures. This tool is used in attacking a
variety of singular convergence, or approximation, problems. We will adopt
this viewpoint in our study of wavelet approximations. Traditionally, the
setting for wavelet questions has included assumptions concerning continuity,
or some kind of differentiability. In contrast, we shall work almost entirely
in the measurable category. One advantage of our present approach is that we
stay in the measurable category when addressing problems from multiresolution
analysis (MRA). Earlier work on the use of probability in wavelets includes
that of R.F. Gundy et al.; see \cite{DGH00,DuJo05b,Fal92,Gun66,Gun99,Gun00,Gun04,GuKa00}.

Our present viewpoint is more general than \cite{DGH00}, and it starts with
the random walks naturally associated with a measure space $\left(
X,\mathcal{B}\right)  $ and a given measurable onto map $\sigma\colon
X\rightarrow X$ such that $\#\sigma^{-1}\left(  \left\{  x\right\}  \right)
=N$ for all $x\in X$, where $N$, $2\leq N<\infty$, is fixed. Iteration of the
branches
\begin{equation}
\sigma^{-1}\left(  \left\{  x\right\}  \right)  :=\left\{  \,y\in X\mid
\sigma\left(  y\right)  =x\,\right\}  \label{eqInt.1}%
\end{equation}
then yields a combinatorial tree. If
\[
\omega=\left(  \omega_{1},\omega_{2},\dots\right)  \in\Omega:=\left\{
0,1,\dots,N-1\right\}  ^{\mathbb{N}},
\]
an associated path may be thought of as an infinite extension of the finite
walks
\[
\tau_{\omega_{n}}\cdots\tau_{\omega_{2}}\tau_{\omega_{1}}x
\]
starting at $x$, where $\left(  \tau_{i}\right)  $, $i=0,1,\dots,N-1$, is a
system of inverses, i.e., where
\begin{equation}
\sigma\circ\tau_{i}=\operatorname*{id}\nolimits_{X},\qquad0\leq i<N.
\label{eqInt.2}%
\end{equation}
If $W\colon X\rightarrow\left[  \,0,1\,\right]  $ is a given measurable
function such that%
\begin{equation}
\sum_{\makebox[0pt]{\hss$\scriptstyle y\colon\sigma\left(
y\right)  =x$\hss}}\;W\left(  y\right)  =1, \label{eqInt.3}%
\end{equation}
then an associated measure $P_{x}$ on $\Omega$ may be defined as follows.
Suppose some function $f\in C\left(  \Omega\right)  $ depends only on a finite
number of coordinates, say $\omega_{1},\dots,\omega_{n}$, then set
\begin{equation}
\int_{\Omega}f\,dP_{x}=\;\;\;\;\sum
_{\makebox[0pt]{\hss$\scriptstyle\left(  \omega_{1},\dots,\omega_{n}\right)$\hss}}%
\;f\left(  \omega_{1},\dots,\omega_{n}\right)  W\left(  \tau_{\omega_{1}%
}x\right)  W\left(  \tau_{\omega_{2}}\tau_{\omega_{1}}x\right)  \cdots
W\left(  \tau_{\omega_{n}}\cdots\tau_{\omega_{1}}x\right)  .\!\!\!\!
\label{eqInt.4}%
\end{equation}
Extensions of this formula to $\Omega$ can be done in a number of ways: see
the cited references.

A special feature of this construction, which will be explored in the present
monograph, is that of attractive convergence properties for infinite products
of the form%
\begin{equation}
\prod_{\makebox[0pt]{\hss$\scriptstyle n,\omega _{1},\dots ,\omega _{n}$\hss}}%
\;W\left(  \tau_{\omega_{1}}x\right)  \cdots W\left(  \tau_{\omega_{n}}%
\cdots\tau_{\omega_{1}}x\right)  \label{eqInt.5}%
\end{equation}
over certain subsets of $\Omega$. As it turns out, these infinite products are
determined by the measures $\left(  P_{x}\right)  _{x\in X}$, and by the
\emph{Ruelle transition operator}%
\begin{equation}
\left(  R_{W}g\right)  \left(  x\right)  :=\;\;\;\;\sum
_{\makebox[0pt]{\hss$\scriptstyle y\colon\sigma\left(
y\right)  =x$\hss}}\;W\left(  y\right)  g\left(  y\right)  ,\qquad g\in
L^{\infty}\left(  X\right)  . \label{eqInt.6}%
\end{equation}

The operator $R$ in (\ref{eqInt.6}) is called the transition operator, the
Ruelle operator, or the Ruelle--Perron--Frobenius operator, and it will play a
major role in what follows.

Many problems in dynamics are governed by transition probabilities $W$, and
$P_{x}$, and by an associated transition operator $R_{W}$ as in (\ref{eqInt.6}%
). Wavelet theory is a case in point, and we show that fundamental convergence
questions for wavelets, and properties of the solutions to (\ref{eqInt.7}),
depend on the positive solutions $h$ to the eigenvalue problem $R_{W} (h) =
h$.
In Theorem
\ref{CorIntSam.4}, we show that there is a particular solution $h$, see (\ref{eqInt.14}), which
determines the issue of pointwise convergence of the infinite product.
We refer to
the discussion around
(\ref{eqInt.14})--(\ref{eqInt.15}) below. The solutions $h$
are called \emph{harmonic}, and the function $h$ in (\ref{eqInt.14}) is a
special harmonic function
which will play a
central role in our main result, Theorem \ref{CorIntSam.4} below.

\section{\label{IntMul}Multiresolutions}

The multiresolution approach to wavelets involves functions on $\mathbb{R}$.
It begins with the fixed-point problem%
\begin{equation}
\varphi\left(  t\right)  =N\sum_{k\in\mathbb{Z}}a_{k}\varphi\left(
Nt-k\right)  ,\qquad t\in\mathbb{R}, \label{eqInt.7bis}%
\end{equation}
where a given sequence $\left(  a_{k}\right)  _{k\in\mathbb{Z}}$ is chosen
with special filtering properties, e.g., quadrature mirror filters; see
\cite{DuJo03,Gro01,HeWe96,JMR01,Jor01a,Mal98,MeCo97,Wal02}. The equation (\ref{eqInt.7bis}) is called the
scaling identity, and the $a_{k}$'s the response coefficients, or the masking
coefficients. Introducing the Fourier series%
\begin{equation}
m\left(  x\right)  =\sum_{k\in\mathbb{Z}}a_{k}e^{-i2\pi kx} \label{eqInt.8}%
\end{equation}
and Fourier transform%
\begin{equation}
\hat{\varphi}\left(  x\right)  =\int_{\mathbb{R}}e^{-i2\pi tx}\varphi\left(
t\right)  \,dt, \label{eqInt.9}%
\end{equation}
we get the relation%
\begin{equation}
\hat{\varphi}\left(  x\right)  =m\left(  \frac{x}{N}\right)  \hat{\varphi
}\left(  \frac{x}{N}\right)  ,\qquad x\in\mathbb{R}, \label{eqInt.10per}%
\end{equation}
which suggests a closer inspection of the infinite products%
\begin{equation}
\prod_{n=1}^{\infty}m\left(  \frac{x}{N^{n}}\right)  . \label{eqInt.10bis}%
\end{equation}
Since we shall want solutions $\varphi$ to (\ref{eqInt.7bis}) which are in
$L^{2}\left(  \mathbb{R}\right)  $, (\ref{eqInt.10per})--(\ref{eqInt.10bis})
suggest the corresponding convergence questions for the function
$W:=\left\vert m\right\vert ^{2}$.

When $\left(  P_{x}\right)  _{x\in\left[  \, 0,1\,\right]  }$ is the family of
measures on $\Omega$ corresponding to $W=\left\vert m\right\vert ^{2}$, then
the formal infinite product%
\begin{equation}
\left\vert \hat{\varphi}\left(  x\right)  \right\vert ^{2}=\prod_{n=1}%
^{\infty}W\left(  \frac{x}{N^{n}}\right)\label{eqInt.10ter}%
\end{equation}
is $P_{x}\left(  \left\{  \left(  0,0,0,\dots\right)  \right\}  \right)  $,
i.e., the measure of the singleton $\left(  0,0,\dots\right)  $ (an infinite
string of zeroes) in $\left\{  0,\dots,N-1\right\}  ^{\mathbb{N}}$.
Our main theorem (Theorem \ref{CorIntSam.4}) gives a necessary and sufficient
condition for the pointwise convergence of (\ref{eqInt.10ter}) in a rather general
context. We note that because of assumption (\ref{eqInt.3}), $\varphi$ will automatically be
in $L^2(\mathbb{R})$, once pointwise convergence is established.
There is a
natural way (based on Euclid's algorithm%
\index{algorithm!Euclid's}%
) of embedding $\mathbb{Z}$ into
\begin{equation}
\Omega=\left\{  0,\dots,N-1\right\}  ^{\mathbb{N}}\times\left\{
0,\dots,N-1\right\}  ^{\mathbb{N}} \label{eqInt.11}%
\end{equation}
such that%
\begin{equation}
P_{x}\left(  \mathbb{Z}\right)  =\sum_{k\in\mathbb{Z}}\left\vert \hat{\varphi
}\left(  x+k\right)  \right\vert ^{2}. \label{eqInt.12}%
\end{equation}

Even though the measures $P_x$ (Section \ref{SomDef}) are defined \emph{a priori} on
the over-countable probability space, the surprise is that they are in fact
supported only on a fixed thin (countable) subset of $\Omega$.

\begin{remark}
\label{RemInt.1}Note that, in general, it is not at all clear that the
measures $\left(  P_{x}\right)  $, $x\in X$, on $\Omega$ should even have
atoms. Typically, they don't! But if atoms exist, i.e., when there are points
$\omega\in\Omega$ such that $P_{x}\left(  \left\{  \omega\right\}  \right)
>0$, we note that this yields convergence of an associated infinite product.
Let $\mathbb{N}_{0}:=\left\{  0,1,2,\dots\right\}  =\mathbb{N}\cup\left\{
0\right\}  $. Using Euclid%
\index{algorithm!Euclid's}%
, and the $N$-adic expansion%
\begin{equation}
k=i_{1}+i_{2}N+\dots+i_{n}N^{n-1}\text{\qquad for }k\in\mathbb{N}_{0},
\label{eqInt.12half}%
\end{equation}
we see that the points
\[
\omega\left(  k\right)  =(\;i_{1},\;\dots,\;i_{n},\;\underbrace
{0,\;0,\;0,\;\dots}_{\!\!\!\infty\text{ string of zeroes}\!\!\!}\;)
\]
represent a copy of $\mathbb{N}_{0}$ sitting in $\Omega$. With the
identification $k\leftrightarrow\omega\left(  k\right)  $, we set%
\begin{equation}
P_{x}\left(  \mathbb{N}_{0}\right)  =\sum_{k=0}^{\infty}P_{x}\left(  \left\{
\omega\left(  k\right)  \right\}  \right)  . \label{eqInt.13}%
\end{equation}
But in general, this function $P_{x}\left(  \mathbb{N}_{0}\right)  $ might be
zero. Our first observation about
\begin{equation}
h\left(  x\right)  :=P_{x}\left(  \mathbb{N}_{0}\right)  ,\qquad x\in X,
\label{eqInt.14}%
\end{equation}
is that it solves the eigenvalue problem%
\begin{equation}
R_{W}h=h. \label{eqInt.15}%
\end{equation}
We say that $h$ is a minimal harmonic function relative to $R_{W}$. 
Note that in general $h=0$ may happen!
\end{remark}

\begin{remark}
\label{RemInt.2}Let $X$, $\mathcal{B}$, $\sigma$, $\tau_{0}$, $\dots$,
$\tau_{N-1}$, and $W$ be as described above, and let $\left(  P_{x}\right)
_{x\in X}$ be the corresponding transition probabilities. Let
\[
\mathbf{0}= \underset{\!\infty\text{ string of zeroes}\!}{(\;0,\;0,\;0,\;\dots
\;)} \in\Omega=\left\{  0,1,\dots,N-1\right\}  ^{\mathbb{N}}.
\]
While in general, often $P_{x}\left(  \left\{  \mathbf{0}\right\}  \right)
=0$, the case when $P_{x}\left(  \left\{  \mathbf{0}\right\}  \right)  >0$ is
important. The condition $P_{x}\left(  \left\{  \mathbf{0}\right\}  \right)
>0$ is a way of making precise sense of the infinite product%
\begin{equation}
P_{x}\left(  \left\{  \mathbf{0}\right\}  \right)  =\prod_{n=1}^{\infty
}W\left(  \tau_{0}^{n}x\right)  . \label{eqInt.16}%
\end{equation}
If, for example, $\lim_{n\rightarrow\infty}W\left(  \tau_{0}^{n}x\right)  <1$,
then it is immediate from \textup{(\ref{eqInt.16})} that $P_{x}\left(
\left\{  \mathbf{0}\right\}  \right)  =0$.

Suppose $P_{x}\left(  \left\{  \mathbf{0}\right\}  \right)  >0$. Then it
follows that%
\begin{equation}
P_{x}\left(  \left\{  \left(  i_{1},\dots,i_{n},\mathbf{0}\right)  \right\}
\right)  =W\left(  \tau_{i_{1}}x\right)  \cdots W\left(  \tau_{i_{n}}%
\cdots\tau_{i_{1}}x\right)  P_{\tau_{i_{n}}\cdots\tau_{i_{1}}x}\left(
\left\{  \mathbf{0}\right\}  \right)  .\!\! \label{eqInt.17}%
\end{equation}
Using \textup{(\ref{eqInt.12half})}, we shall identify $k\in\mathbb{N}_{0}$
with the point $\omega\left(  k\right)  \in\Omega$, and write $P_{x}\left(
\left\{  k\right\}  \right)  $ for the expression in \textup{(\ref{eqInt.17})}.
\end{remark}

An important question for dyadic wavelets in $L^{2}\left(  \mathbb{R}\right)
$ is the issue of when these wavelets form \emph{orthonormal bases} (ONB's). A
dyadic wavelet function $\psi\in L^{2}\left(  \mathbb{R}\right)  $ generates
an ONB if the double-indexed family
\begin{equation}
\left\{  \,2^{n/2}\psi\left(  2^{n}t-k\right)  \bigm| n,k\in\mathbb{Z}%
\,\right\}  \label{eqIntMul.16}%
\end{equation}
satisfies (\ref{ONB(1)}) and (\ref{ONB(2)}) below:

\begin{enumerate}
\renewcommand{\theenumi}{\roman{enumi}}

\item \label{ONB(1)}$\displaystyle\int_{\mathbb{R}}\overline{\psi_{n,k}\left(
t\right)  }\,\psi_{m,l}\left(  t\right)  \,dt=\delta_{n,m}\delta_{k,l}$, with%
\begin{equation}
\psi_{n,k}\left(  t\right)  :=2^{n/2}\psi\left(  2^{n}t-k\right)  ,
\label{eqIntMul.17}%
\end{equation}

\end{enumerate}

\noindent\setcounter{saveenumi}{\value{enumi}}and

\begin{enumerate}
\renewcommand{\theenumi}{\roman{enumi}}\setcounter{enumi}{\value{saveenumi}}

\item \label{ONB(2)}the closed linear span of $\left\{  \,\psi_{n,k}\mid
n,k\in\mathbb{Z}\,\right\}  $ is $L^{2}\left(  \mathbb{R}\right)  $.
\end{enumerate}

In our analysis of the scaling identity%
\begin{equation}
\varphi\left(  t\right)  =2\sum_{k\in\mathbb{Z}}a_{k}\varphi\left(
2t-k\right)  \label{eqIntMul.18}%
\end{equation}
(a special case of (\ref{eqInt.7bis})), we will be looking at two functions
$\varphi$ and $\psi$; the second one may be taken to be%
\begin{equation}
\psi\left(  t\right)  =2\sum_{k\in\mathbb{Z}}\left(  -1\right)  ^{k+1}\bar
{a}_{1-k}\varphi\left(  2t-k\right)  . \label{eqIntMul.19}%
\end{equation}
This analysis is the approach to wavelets which goes under the name of
\emph{multiresolution analysis}. The function $\psi$ which is used in
(\ref{eqIntMul.16})--(\ref{eqIntMul.17}) is the solution to (\ref{eqIntMul.19}%
). The two standing conditions which are placed on the numbers $\left(
a_{k}\right)  _{k\in\mathbb{Z}}$, called \emph{masking coefficients}, are%
\begin{equation}
\sum_{k\in\mathbb{Z}}\bar{a}_{k}a_{k+2n}=\frac{1}{2}\delta_{0,n},\qquad
n\in\mathbb{Z}, \label{eqIntMul.20}%
\end{equation}
and%
\begin{equation}
\sum_{k\in\mathbb{Z}}a_{k}=1. \label{eqIntMul.21}%
\end{equation}
These conditions in themselves do not imply orthonormality in
(\ref{eqIntMul.17}), but only the following much weaker property:%
\begin{equation}
\mathop{\sum\!\sum}_{n,k\in\mathbb{Z}}\left\vert \left\langle \,\psi_{n,k}\mid
f\,\right\rangle \right\vert ^{2}=\left\Vert f\right\Vert ^{2}=\int
_{\mathbb{R}}\left\vert f\left(  t\right)  \right\vert ^{2}\,dt,\qquad f\in
L^{2}\left(  \mathbb{R}\right)  . \label{eqIntMul.22}%
\end{equation}
A system of functions $\left(  \psi_{n,k}\right)  $ satisfying
(\ref{eqIntMul.22}) is called a \emph{Parseval frame}, or a \emph{normalized
tight frame}.

Given (\ref{eqIntMul.20})--(\ref{eqIntMul.21}), it turns out
that the ONB property for the wavelet is equivalent to either one
of the following two conditions for the normalized scaling function $\varphi
$\kern0.5pt:%
\begin{equation}
\int_{\mathbb{R}}\overline{\varphi\left(  t\right)  }\,\varphi\left(
t-k\right)  \,dt=\delta_{0,k}, \label{eqIntMul.23}%
\end{equation}
or%
\begin{equation}
\left\Vert \varphi\right\Vert _{L^{2}\left(  \mathbb{R}\right)  }=1.
\label{eqIntMul.24}%
\end{equation}

\section{\label{IntSam}Sampling}

In this section and the next, we study an intriguing relationship between the
following three problems:

\begin{enumerate}
\renewcommand{\theenumi}{\arabic{enumi}}

\item \label{IntSam1} When does the scaling identity (\ref{eqInt.7bis}) have
$L^{2}$ solutions?

\item \label{IntSam2} How may the transition probabilities $P_{x}$ be used in
\emph{sampling} certain functions at the points $x + k$ as $k$ runs over a set
of integers?

\item \label{IntSam3} When is the infinite product (\ref{eqInt.10bis})
pointwise convergent?
\end{enumerate}

In the wavelet applications, $X=\left[  \, 0,1\,\right]  $, and the system
$\sigma,\tau_{0},\dots,\tau_{N-1}$ is as follows:%
\begin{equation}
\left\{
\begin{aligned} &\sigma\left( x\right) =Nx\bmod1,\\ &\tau_{j}\left( x\right) =\frac{x+j}{N},\qquad j=0,1,\dots,N-1. \end{aligned} \right.
\label{eqInt.18}%
\end{equation}

\begin{lemma}
\label{LemInt.3}Setting $F\left(  x\right)  :=P_{x}\left(  \left\{
\mathbf{0}\right\}  \right)  $, and%
\begin{equation}
k=i_{1}+i_{2}N+\dots+i_{n}N^{n-1}\qquad(\in\mathbb{N}_{0}), \label{eqInt.19}%
\end{equation}
we get the formula%
\begin{equation}
P_{x}\left(  \left\{  k\right\}  \right)  =F\left(  x+k\right)  ,
\label{eqInt.20}%
\end{equation}
where we have identified $k$ with the point $\omega\left(  k\right)  :=\left(
i_{1},\dots,i_{n},\mathbf{0}\right)  $ in $\Omega$.
\end{lemma}

\begin{proof}
To see this, identify functions on $\left[  \, 0,1\,\right]  $ with
$1$-periodic functions on $\mathbb{R}$, and note that the second formula in
(\ref{eqInt.18}) yields $\tau_{i_{n}}\cdots\tau_{i_{1}}\left(  x\right)
=\left(  x+k\right)  /N^{n}$ where $k$ is given by (\ref{eqInt.19}). Hence, if
$1\leq s\leq n$, then%
\[
W\left(  \tau_{i_{s}}\cdots\tau_{i_{1}}\left(  x\right)  \right)  =W\left(
\frac{x+i_{1}N+\cdots+i_{s}N^{s-1}}{N^{s}}\right)  =W\left(  \frac{x+k}{N^{s}%
}\right)  .
\]
It follows that (\ref{eqInt.20}) is really just a rewrite of (\ref{eqInt.17}).
The right-hand side of (\ref{eqInt.17}) yields%
\begin{equation}
W\left(  \frac{x+k}{N}\right)  \cdots W\left(  \frac{x+k}{N^{n}}\right)
F\left(  \frac{x+k}{N^{n}}\right)  =F\left(  x+k\right)  ,
\label{eqInt.20half}%
\end{equation}
which is the desired conclusion.
\end{proof}

\begin{remark}
\label{RemInt.4}To extend $P_{x}\left(  \,\cdot\,\right)  $ from
$\mathbb{N}_{0}$ to $\mathbb{Z}$, recall that $k\in\mathbb{N}_{0}$ is
identified with the singleton $\omega\left(  k\right)  =\left(  i_{1}%
,\dots,i_{n},\mathbf{0}\right)  $ via \textup{(\ref{eqInt.19})}. Now, if
$-N^{n}\leq k<0$, then set
\begin{equation}
P_{x}\left(  \left\{  k\right\}  \right)  :=P_{x}\left(  \left\{
\omega\left(  N^{n+1}+k\right)  \right\}  \right)  . \label{eqInt.21}%
\end{equation}

To help the reader gain some intuitive feeling for the conclusion in Lemma
\textup{\ref{LemInt.3}}, observe that the right-hand side of
\textup{(\ref{eqInt.20})} represents a \emph{sampling} of the function $F$ at
the integral translates on $\mathbb{R}$, starting at $x$, i.e., $x+k$.
Obviously, different subsets of $\Omega$ would yield different sets of
sampling points for $F$, including nonuniformly distributed sampling points;
see \cite{AlGr01}.

Starting with Shannon \cite{Sha49}, the theory of sampling has emerged as a
significant tool in signal processing; see, e.g., the beautifully written
survey \cite{AlGr01} as well as the references cited therein. Thus, in a
general context, our formula \textup{(\ref{eqInt.20})} offers a probabilistic
prescription for sampling of functions on the real line, and at the same time
it stresses the `random' feature of sampling.
\end{remark}

\section{\label{ConThm}A convergence theorem for infinite products}

We now show how this viewpoint from sampling theory is closely related to some
fundamental properties of the measures $P_{x}$. In particular, our Theorem
\ref{CorIntSam.4} gives a necessary and sufficient condition for pointwise
convergence of the infinite product (\ref{eqInt.10bis}), or more generally
(\ref{eqInt.5}), with the condition for convergence stated in terms of the
harmonic function $h$ of (\ref{eqInt.14}). The relationship between $h$ and
the measure family $P_{x}$ is studied more systematically below.

Let $A\subset\mathbb{Z}$. Returning to (\ref{eqInt.20}), we set%
\begin{equation}
P_{x}\left(  A\right)  :=\sum_{k\in A}P_{x}\left(  \left\{  k\right\}
\right)  =\sum_{k\in A}F\left(  x+k\right)  . \label{eqIntSam.5}%
\end{equation}
As in Lemma \ref{LemInt.3}, the number $N$, $N\geq2$, and the function $W$ are
given. The measures $P_{x}$ are constructed from these data using
(\ref{eqInt.4}), and we have the two functions%
\begin{equation}
F\left(  x\right)  :=P_{x}(\;\{\;(\;\underbrace{0,\;0,\;0,\;\dots
}_{\!\!\!\infty\text{ string of zeroes}\!\!\!}\;)\;\}\;) \label{eqIntSam.6}%
\end{equation}
and%
\begin{equation}
h\left(  x\right)  :=P_{x}\left(  \mathbb{Z}\right)  ,\qquad x\in\mathbb{R}.
\label{eqIntSam.7}%
\end{equation}
Finally, for $k\in\mathbb{N}$, set%
\begin{equation}
N^{k}\mathbb{Z}:=\left\{  \,N^{k}j\mid j\in\mathbb{Z}\,\right\}  .
\label{eqIntSam.8}%
\end{equation}
Using (\ref{eqInt.21}), we see that $N^{k}\mathbb{Z}$ is represented in
$\Omega=\left\{  0,1,\dots,N-1\right\}  ^{\mathbb{N}}$ as%
\begin{equation}
(\;\underbrace{0,\;\dots,\;0}_{k\text{ zeroes}},\;\underbrace{\omega
_{1},\;\omega_{2},\;\dots,\;.\,}_{%
\begin{smallmatrix}
\text{a finite string}\\
\text{of symbols}\\
\omega_{i}\in\left\{  0,\dots,N-1\right\}
\end{smallmatrix}
},\;\underbrace{0,\;0,\;0,\;\dots}_{\!\!\!\infty\text{ string of
zeroes}\!\!\!}\;). \label{eqIntSam.9}%
\end{equation}
An infinite string of zeroes will be denoted $\mathbf{0}$.

\begin{lemma}
\label{CorIntSam.3}Let $N$, $W$, $F$, $P_{x}$, and $h$ be as described above;
see \textup{(\ref{eqIntSam.6})--(\ref{eqIntSam.7})}. Then $h$ satisfies the
following cocycle identity:%
\begin{equation}
h\left(  x\right)  W\left(  x\right)  =P_{Nx}\left(  N\mathbb{Z}\right)
,\qquad x\in\mathbb{R}. \label{eqIntSam.10}%
\end{equation}
In particular, if $h\left(  x\right)  \equiv1$ $\mathrm{a.e.}\,x\in\mathbb{R}%
$, then we recover the function $W$ from the transition probabilities $P_{x}$,
$\mathrm{a.e.}\,x\in\mathbb{R}$.
\end{lemma}

\begin{proof}
We calculate the left-hand side in (\ref{eqIntSam.10}), using the earlier
equations:%
\begin{align*}
h\left(  x\right)  W\left(  x\right)   &  \underset
{\textstyle\text{\makebox[\byperiodicity]{$\begin{smallmatrix}  \text{by (\ref{eqInt.20})}\\
\text{and (\ref{eqIntSam.7})}\end{smallmatrix}$}}}{=}W\left(  x\right)
\sum_{j\in\mathbb{Z}}F\left(  x+j\right) \\
&  \underset{\text{by periodicity}}{=}\sum_{j\in\mathbb{Z}}W\left(
x+j\right)  F\left(  x+j\right) \\
&  \underset{\text{\makebox[\byperiodicity]{by (\ref{eqInt.20half})}}}{=}%
\sum_{j\in\mathbb{Z}}F\left(  N\left(  x+j\right)  \right) \\
&  \underset{\text{\makebox[\byperiodicity]{~}}}{=}\sum_{j\in\mathbb{Z}%
}F\left(  Nx+Nj\right) \\
&  \underset{\text{\makebox[\byperiodicity]{by (\ref{eqInt.20})}}}{=}%
\sum_{j\in\mathbb{Z}}P_{Nx}\left(  \left\{  Nj\right\}  \right) \\
&  \underset{\text{\makebox[\byperiodicity]{by (\ref{eqIntSam.5})}}}{=}%
P_{Nx}\left(  N\mathbb{Z}\right)  .
\end{align*}
This is the desired identity (\ref{eqIntSam.10}), and the proof is completed.
\end{proof}

\begin{theorem}
\label{CorIntSam.4}Let $N$, $W$, $F$, $P_{x}$, and $h$ be as described above.
Let $x\in\mathbb{R}$, and suppose that 
$P_{x}\left(  \left\{  \mathbf{0}\right\}  \right)  >0$. Then the
following two conditions are equivalent.

\begin{enumerate}
\renewcommand{\theenumi}{\alph{enumi}}

\item \label{CorIntSam.4(1)}The limit on the right-hand side below exists, and%
\[
F\left(  x\right)  =\lim_{n\rightarrow\infty}\prod_{k=1}^{n}W\left(  \frac
{x}{N^{k}}\right)  .
\]

\item \label{CorIntSam.4(2)}The limit on the left-hand side below exists, and%
\[
\lim_{n\rightarrow\infty}h\left(  \frac{x}{N^{n}}\right)  =1.
\]

\end{enumerate}
\end{theorem}

\begin{proof}
(\ref{CorIntSam.4(1)})${}\Rightarrow{}$(\ref{CorIntSam.4(2)}). An iteration of
the identity (\ref{eqIntSam.10}) in Lemma \ref{CorIntSam.3} above yields%
\begin{equation}
P_{x}\left(  N^{k}\mathbb{Z}\right)  =\left(  \prod_{j=1}^{k}W\left(  \frac
{x}{N^{j}}\right)  \right)  h\left(  \frac{x}{N^{k}}\right)  .
\label{eqIntSam.11}%
\end{equation}
Using (\ref{eqIntSam.9}), and working in $\Omega$, we find%
\begin{equation}
\bigcap_{k\in\mathbb{N}}N^{k}\mathbb{Z}=\left\{  \mathbf{0}\right\}  .
\label{eqIntSam.12}%
\end{equation}
An application of a standard result in measure theory \cite[Theorem 1.19(e),
p.~16]{Rud87} now yields existence of the following limit:%
\begin{equation}
\lim_{k\rightarrow\infty}P_{x}\left(  N^{k}\mathbb{Z}\right)  =P_{x}\left(
\left\{  \mathbf{0}\right\}  \right)  \underset{\text{by (\ref{eqIntSam.6})}%
}{=}F\left(  x\right)  . \label{eqIntSam.13}%
\end{equation}
Since (\ref{CorIntSam.4(1)}) is assumed, and $F\left(  x\right)  >0$, we
conclude that the limit $h\left(  x/N^{k}\right)  $, for $k\rightarrow\infty$,
must exist as well, and further that%
\[
F\left(  x\right)  =F\left(  x\right)  \lim_{k\rightarrow\infty}h\left(
\frac{x}{N^{k}}\right)  .
\]
Using again $F\left(  x\right)  >0$, we finally conclude that
(\ref{CorIntSam.4(2)}) holds.

(\ref{CorIntSam.4(2)})${}\Rightarrow{}$(\ref{CorIntSam.4(1)}). Recall that
formula (\ref{eqIntSam.11}) holds in general. If (\ref{CorIntSam.4(2)}) is
assumed, we then conclude that the limit $\prod_{j=1}^{k}W\left(
x/N^{j}\right)  $ exists as $k\rightarrow\infty$. The limit on the left-hand
side in (\ref{eqIntSam.11}) exists and is $F\left(  x\right)  $. As a result,
we get%
\[
F\left(  x\right)  =\lim_{k\rightarrow\infty}\prod_{j=1}^{k}W\left(  \frac
{x}{N^{j}}\right)  \lim_{k\rightarrow\infty}h\left(  \frac{x}{N^{k}}\right)
\;\underset{\text{using (\ref{CorIntSam.4(2)})}}{=}\;\lim_{k\rightarrow\infty
}\prod_{j=1}^{k}W\left(  \frac{x}{N^{j}}\right)  ,
\]
which is (\ref{CorIntSam.4(1)}).
\end{proof}

\section{\label{SomDef}Ruelle's wavelet transition operator}

The
next definitions (Definitions \ref{DefTra.1}) and the lemma (Lemma \ref{LemTra.2}) give us the
precise details behind the two critical notions from Theorem \ref{CorIntSam.4} above;
i.e., the cocycles, and the random walk measures $P_x$. An existence and
uniqueness theorem for $P_x$ is then presented in Section \ref{KolCon}.

\begin{definitions}
\label{DefTra.1}\renewcommand{\theenumi}{\alph{enumi}}\ 

\begin{enumerate}
\item \label{DefTra.1(1)}The \emph{Ruelle operator} $R=R_{W}$ is defined by%
\begin{equation}
\left(  Rf\right)  \left(  x\right)  =\sum_{y\in X,\,\sigma\left(  y\right)
=x}\!\!\!\!W\left(  y\right)  f\left(  y\right)  ,\qquad x\in X,\;f\in
L^{\infty}\left(  X\right)  , \label{eqTra.4}%
\end{equation}
and maps $L^{\infty}\left(  X\right)  $ into itself.

\item \label{DefTra.1(2)}Let $\Omega$ be the compact Cartesian product
\begin{equation}
\Omega=\mathbb{Z}_{N}^{\mkern4mu\mathbb{N}}=\left\{  0,\dots,N-1\right\}
^{\mathbb{N}}=\prod_{1}^{\infty}\left\{  0,\dots,N-1\right\}  .
\label{eqTra.5}%
\end{equation}

\item \label{DefTra.1(3)}A bounded measurable function $V\colon X\times
\Omega\rightarrow\mathbb{C}$ is said to be a \emph{cocycle} if%
\begin{equation}
V\left(  x,\left(  \omega_{1},\omega_{2},\dots\right)  \right)  =V\left(
\tau_{\omega_{1}}\left(  x\right)  ,\left(  \omega_{2},\omega_{3}%
,\dots\right)  \right)  \label{eqTra.6}%
\end{equation}
for all $\omega=\left(  \omega_{1},\omega_{2},\dots\right)  \in\Omega$.

\item \label{DefTra.1(4)}A function $h\colon X\rightarrow\mathbb{C}$ is said
to be \emph{harmonic}, or $R_{W}$\emph{-harmonic}, if%
\begin{equation}
R_{W}h=h. \label{eqTra.7}%
\end{equation}

\item \label{DefTra.1(5)}Let $n\in\mathbb{N}$, and let $i_{1},\dots,i_{n}%
\in\mathbb{Z}_{N}$. Then the subset%
\begin{equation}
A\left(  i_{1},\dots,i_{n}\right)  :=\left\{  \,w\in\Omega\mid\omega_{1}%
=i_{1},\;\dots,\;\omega_{n}=i_{n}\,\right\}  \label{eqTra.8}%
\end{equation}
is called a \emph{cylinder set}.
\end{enumerate}
\end{definitions}

We shall use the following correspondence (see \cite{Jor05} and \cite{DuJo05b}
for details and proofs) between the cocycles $V$ from (\ref{DefTra.1(3)}) and the harmonic
functions $h$ from (\ref{DefTra.1(4)}): If $V$ is given as in (\ref{DefTra.1(3)}), then the function $h(x):= P_x[
V(x,\,\cdot\,)]$ has the properties in (\ref{DefTra.1(4)}). Conversely, for every $h$ satisfying (\ref{DefTra.1(4)}),
including (\ref{eqTra.7}), there is a martingale limit which lets us recover a cocycle
$V$, $P_x$ a.e., such that $h(x):= P_x[ V(x,\,\cdot\,)]$. In the present discussion, we
are interested in the cocycles $V$ which arise as indicator functions  $\chi_S^{{}}$
for certain cyclic and invariant subsets $S$ of $\Omega$. We refer to Section \ref{KolCon}
below.

\subsection{\label{ExiMea}Existence of the measures $P_{x}$}

The cylinder sets generate the topology of $\Omega$, and its Borel
sigma-algebra. In determining Radon measures on $\Omega$, it is therefore
convenient to first specify them on cylinder sets. This approach was initiated
by Kolmogorov \cite{Kol77}; see also Nelson \cite{Nel69}. Recall that $\Omega$
is compact in the Tychonoff topology, and that we may use the
Stone--Weierstra\ss \ theorem on $C\left(  \Omega\right)  ={}$the algebra of
all continuous functions on $\Omega$.

\begin{lemma}
\label{LemTra.2}Let $X$, $\mathcal{B}$, $\mu$, $\sigma$, $\tau_{0}$, $\dots$,
$\tau_{N-1}$, and $W$ be given as described above. We make the following more
restrictive assumption on $W$:%
\begin{equation}
\sum_{\makebox[0pt]{\hss$\scriptstyle y\in X,\,\sigma\left(
y\right)  =x$\hss}}\;W\left(  y\right)  =1\qquad\mathrm{a.e.}\;x\in X.
\label{eqTra.9}%
\end{equation}
Then for every $x\in X$ there is a unique positive Radon probability measure
$P_{x}$ on $\Omega$ such that%
\begin{equation}
P_{x}\left(  A\left(  i_{1},\dots,i_{n}\right)  \right)  =W\left(  \tau
_{i_{1}}x\right)  W\left(  \tau_{i_{2}}\tau_{i_{1}}x\right)  \cdots W\left(
\tau_{i_{n}}\cdots\tau_{i_{1}}x\right)  . \label{eqTra.10}%
\end{equation}
\end{lemma}

The main fact about the Ruelle operator
\[
R_{W}f\left( x\right) =\sum_{i}W\left( \tau_{i}x\right) f\left( \tau_{i}x\right) 
\]
is this \cite{DuJo05b}:
Under suitable conditions on $W$,
there is a unique probability measure (Ruelle measure $\nu$) satisfying
\[
\nu \cdot R_{W}=\nu,
\]
and a unique continuous minimal eigenfunction $h_{\min}$,
\[
R_{W}h_{\min}=h_{\min}.
\]
This function $h_{\min}$ is minimal in the ordered convex set
\[
\left\{\,h\in C\left( X\right) \Bigm| 0\leq h\leq 1,\;R_{W}h=h,\;\int_{X}h\,d\nu=1\right\},
\]
and moreover
\[
h_{\min}\left( x\right) =P_{x}\left( \mathbb{Z}\right) .
\]

\begin{remark}
\label{RemTra.2bis}It turns out that the general case $\sum_{y}\cdots\leq1$
may be
reduced to \textup{(\ref{eqTra.9})}. So \textup{(\ref{eqTra.9})} is not really
a restriction.
\end{remark}

\begin{proof}
[Proof of Lemma \textup{\ref{LemTra.2}}]If $P$ is a Radon measure on $\Omega$,
we set
\begin{equation}
P\left[  \, f\,\right]  :=\int_{\Omega}f\left(  \omega\right)  \,dP\left(
\omega\right)  \text{\qquad for all }f\in C\left(  \Omega\right)  .
\label{eqTra.10bis}%
\end{equation}
Set%
\begin{multline}
C_{\operatorname*{fin}}\left(  \Omega\right)  =\{\,f\in C\left(
\Omega\right)  \mid\exists\,n\text{ such that }f\left(  \omega\right)
=f\left(  \omega_{1},\dots,\omega_{n}\right)  \text{, }\label{eqTra.11}\\[1pt]
\text{i.e., }f\text{ depends only on the first }n\text{ coordinates in }%
\Omega\,\}.\!\!
\end{multline}
Increasing $n$ in definition (\ref{eqTra.11}), we get an ascending nest of
subalgebras of $C\left(  \Omega\right)  $,%
\begin{equation}
\mathfrak{A}_{1}\subset\mathfrak{A}_{2}\subset\dots\subset\mathfrak{A}%
_{n}\subset\mathfrak{A}_{n+1}\subset\cdots, \label{eqTra.11bis}%
\end{equation}
with%
\[
\overline{\bigcup_{n=1}^{\infty}\mathfrak{A}_{n}}\,=C\left(  \Omega\right)  ,
\]
where $\overline{\,\cdot\,}$ stands for norm-closure. We set%
\[
C_{\operatorname*{fin}}\left(  \Omega\right)  :=\bigcup_{n=1}^{\infty
}\mathfrak{A}_{n}.
\]
An immediate application of Stone--Weierstra\ss \ shows that
$C_{\operatorname*{fin}}\left(  \Omega\right)  $ is uniformly dense in
$C\left(  \Omega\right)  $. Let $x\in X$, and $f\in C_{\operatorname*{fin}%
}\left(  \Omega\right)  $. Suppose
\[
f\left(  \omega\right)  =f\left(  \omega_{1},\dots,\omega_{n}\right)  ,
\]
and set%
\begin{equation}
P_{x}\left[  \, f\,\right]  =\;\;\;\;\;\;\sum
_{\makebox[0pt]{\hss$\scriptstyle\left(  \omega_{1},\dots,\omega_{n}\right)\in\mathbb{Z}_{N}^{\mkern4mu n}$\hss}}%
\;\;W\left(  \tau_{\omega_{1}}x\right)  \cdots W\left(  \tau_{\omega_{n}%
}\cdots\tau_{\omega_{1}}x\right)  f\left(  \omega_{1},\dots,\omega_{n}\right)
. \label{eqTra.12}%
\end{equation}
Note that, if there is some $h\in L^{\infty}\left(  X\right)  $ such that
\[
f\left(  \omega_{1},\dots,\omega_{n}\right)  =h\left(  \tau_{\omega_{n}}%
\cdots\tau_{\omega_{1}}x\right)  ,
\]
then%
\begin{align}
P_{x}\left[  \, f\,\right]   &  =\;\;\;\;\sum
_{\makebox[0pt]{\hss$\scriptstyle\left(  \omega_{1},\dots,\omega_{n}\right)$\hss}}%
\;W\left(  \tau_{\omega_{1}}x\right)  \cdots W\left(  \tau_{\omega_{n}}%
\cdots\tau_{\omega_{1}}x\right)  h\left(  \tau_{\omega_{n}}\cdots\tau
_{\omega_{1}}x\right) \label{eqTra.13}\\
&  =\;\;\;\;\;\;\sum
_{\makebox[0pt]{\hss$\scriptstyle y \in X,\,\sigma^{n}y=x$\hss}}\;\;W\left(
\sigma^{n-1}y\right)  \cdots W\left(  \sigma y\right)  W\left(  y\right)
h\left(  y\right) \nonumber\\
&  =\left(  R_{W}^{n}h\right)  \left(  x\right)  .\nonumber
\end{align}

We now show that $P_{x}\left[  \, f\,\right]  $ is well-defined. This is the
Kolmogorov consistency: we must check that the number $P_{x}\left[  \,
f\,\right]  $ is the same when some $f\in\mathfrak{A}_{n}$ ($\subset
\mathfrak{A}_{n+1}$) is viewed also as an element in $\mathfrak{A}_{n+1}$.
Then
\begin{align*}
f\left(  \omega\right)   &  =f\left(  \omega_{1},\dots,\omega_{n}\right) \\
&  =f\left(  \omega_{1},\dots,\omega_{n},\omega_{n+1}\right)  ,
\end{align*}
and
\begin{align*}
P_{x}\left[  \, f_{n+1}\,\right]   &  =\;\;\;\;\sum
_{\makebox[0pt]{\hss$\scriptstyle  \omega_{1},\dots,\omega_{n+1}$\hss}}%
\;W\left(  \tau_{\omega_{1}}x\right)  \cdots W\left(  \tau_{\omega_{n+1}%
}\cdots\tau_{\omega_{1}}x\right)  f\left(  \omega_{1},\dots,\omega
_{n+1}\right) \\
&  =\;\;\;\sum
_{\makebox[0pt]{\hss$\scriptstyle  \omega_{1},\dots,\omega_{n}$\hss}}%
\;W\left(  \tau_{\omega_{1}}x\right)  \cdots W\left(  \tau_{\omega_{n}}%
\cdots\tau_{\omega_{1}}x\right) \\
&  \qquad\qquad\quad{}\cdot\left(  \;
\vphantom{\sum_{\omega_{n+1}}W\left( \tau_{\omega_{n+1}}\tau_{\omega_{n}}\cdots\tau_{\omega_{1}}x\right) }
\smash{\underbrace{\sum_{\omega_{n+1}}W\left( \tau_{\omega_{n+1}}\tau_{\omega_{n}}\cdots\tau_{\omega_{1}}x\right) }_{=1\text{ by (\ref{eqTra.9})}}}
\; \right)
\vphantom{\underbrace{\sum_{\omega_{n+1}}W\left( \tau_{\omega_{n+1}}\tau_{\omega_{n}}\cdots\tau_{\omega_{1}}x\right) }_{=1\text{ by (\ref{eqTra.9})}}}
f\left(  \omega_{1},\dots,\omega_{n}\right) \\
&  =\;\;\;\sum
_{\makebox[0pt]{\hss$\scriptstyle  \omega_{1},\dots,\omega_{n}$\hss}}%
\;W\left(  \tau_{\omega_{1}}x\right)  \cdots W\left(  \tau_{\omega_{n}}%
\cdots\tau_{\omega_{1}}x\right)  f\left(  \omega_{1},\dots,\omega_{n}\right)
\\
&  =P_{x}\left[  \, f_{n}\,\right]  ,
\end{align*}
as claimed.

The consistency conditions may be stated differently in terms of conditional
probabilities: for $f\in C\left(  \Omega\right)  $, set%
\begin{align}
P_{x}^{\left(  n\right)  }\left[  \,f\,\right]   &  =P_{x}\left[
\,f\mid\mathfrak{A}_{n}\,\right] \label{eqTra.13bis}\\
&  =\;\;\;\;\sum
_{\makebox[0pt]{\hss$\scriptstyle\left(  \omega_{1},\dots,\omega_{n}\right)$\hss}}%
\;W\left(  \tau_{\omega_{1}}x\right)  \cdots W\left(  \tau_{\omega_{n}}%
\cdots\tau_{\omega_{1}}x\right)  f\left(  \omega_{1},\dots,\omega_{n}\right)
.\nonumber
\end{align}
We proved that
\[
P_{x}^{\left(  n\right)  }\left[  \,f\,\right]  =P_{x}^{\left(  n+1\right)
}\left[  \,f\,\right]  \text{\qquad for all }f\in\mathfrak{A}_{n}.
\]
Using now the theorems of Stone--Weierstra\ss \ and Riesz, we get the
existence of the measure $P_{x}$ on $\Omega$. It is clear that it has the
desired properties. In particular, the property (\ref{eqTra.10}) results from
applying (\ref{eqTra.12}) to the function%
\begin{equation}
f\left(  \omega\right)  :=\delta_{i_{1},\omega_{1}}\cdots\delta_{i_{n}%
,\omega_{n}},\qquad\omega\in\Omega, \label{eqTra.14}%
\end{equation}
when the point $\left(  i_{1},\dots,i_{n}\right)  $ is fixed.

These functions, in turn, span a dense subalgebra in $C\left(  \Omega\right)
$ (by Stone--Weier\-stra\ss ), so $P_{x}$ is determined uniquely by
(\ref{eqTra.10}).
\end{proof}

\section{\label{BriOut}Transition measures}

We now compute the Markov transition measures $P_{x}$, as $x$ varies over the
set $X$. The construction of $P_{x}$, and the probability space $\Omega$,
begins with a sequence of measures $P_{x}^{\left(  n\right)  }$, $n =
1,2,\dots$, corresponding to finite paths of length $n$. Then we show that the
infinite path-space measure, Lemma \ref{LemTra.2}, results from an application
of the Kolmogorov extension principle. Our technical analysis involves a
certain transition operator $R$ which generalizes Lawton's wavelet transition
operator; see \cite{Law91a,Law91b}.

We already mentiond how the measures $P_{x}$ serve to prescribe sampling of
functions on the real line, Lemma \ref{LemInt.3}. However, a deeper
understanding of this sampling viewpoint is facilitated by the introduction of
the Perron--Frobenius--Ruelle operator $R$, Definitions \ref{DefTra.1}, and an
associated family of harmonic functions.
For these harmonic functions, we note in Section \ref{SomDef} that there
is a crucial analogue from classical harmonic analysis of Fatou-boundary
function in the present discrete context.  As noted after Definition \ref{DefTra.1}
below, the boundary value functions take the form of certain cocycles, and
the function $h$ from (\ref{eqInt.14}) corresponds to a special $\mathbb{N}_0$-cocycle. The
existence of these cocycles is based on a martingale convergence theorem;
see \cite[Theorem 2.7.1]{Jor05} and \cite{DuJo05b} for additional details.

A fundamental tightness condition for the random walk
model is introduced. The transition probabilities $P_{x}$ live in a universal
probability space $\Omega$, but the essential convergence questions for the
infinite products depend on the the $P_{x}$'s being
supported on a certain copy of $\mathbb{N}_{0}$ (the natural numbers), or of
$\mathbb{Z}$ (the integers). This refers to the natural embedded of
$\mathbb{Z}$ in $\Omega$ which we described above.

\section{\label{KolCon}Kolmogorov's consistency condition}

For general reference, we now make explicit the extension principle of
Kolmogorov \cite{Kol77} in its function-theoretic form.

\begin{lemma}
\label{LemTraNew.pound}\textup{(Kolmogorov)}\enspace Let $N\geq2$ be fixed,
and let
\[
\Omega=\left\{  0,1,\dots,N-1\right\}  ^{\mathbb{N}}.
\]
For $n=1,2,\dots$, let%
\[
P^{\left(  n\right)  }\colon\mathfrak{A}_{n}\rightarrow\mathbb{C}%
\]
be a sequence of linear functionals such that
\textup{(\ref{LemTraNew.pound(1)})--(\ref{LemTraNew.pound(3)})} hold:

\begin{enumerate}
\renewcommand{\theenumi}{\roman{enumi}}

\item \label{LemTraNew.pound(1)}$P^{\left(  n\right)  }\left[  \,
\openone\,\right]  =1$, where $\openone$ denotes the constant function $1$ on
$\Omega$,

\item \label{LemTraNew.pound(2)}$f\in\mathfrak{A}_{n}$, $f\geq0$
pointwise${}\Longrightarrow P^{\left(  n\right)  }\left[  \, f\,\right]
\geq0$,
\end{enumerate}

\noindent\setcounter{saveenumi}{\value{enumi}}and

\begin{enumerate}
\renewcommand{\theenumi}{\roman{enumi}}\setcounter{enumi}{\value{saveenumi}}

\item \label{LemTraNew.pound(3)}$P^{\left(  n\right)  }\left[  \, f\,\right]
=P^{\left(  n+1\right)  }\left[  \, f\,\right]  $ for all $f\in\mathfrak{A}%
_{n}$.
\end{enumerate}

Then there is a unique Borel probability measure $P$ on $\Omega$ such that%
\begin{equation}
P\left[  \, f\,\right]  =P^{\left(  n\right)  }\left[  \, f\,\right]  ,\qquad
f\in\mathfrak{A}_{n}. \label{eqLemTraNew.pound1}%
\end{equation}
Specifically, for $P$, we have the implication%
\begin{equation}
f\in C\left(  \Omega\right)  ,\;f\geq0\text{ pointwise}{}\Longrightarrow
P\left[  \, f\,\right]  \geq0. \label{eqLemTraNew.pound2}%
\end{equation}

\end{lemma}

\begin{remark}
\label{RemTraNew.poundter}Here we have identified positive linear functionals
$P$ on $C\left(  \Omega\right)  $ with the corresponding Radon measures
$\tilde{P}$ on $\Omega$, i.e.,%
\begin{equation}
P\left[  \, f\,\right]  =\int_{\Omega}f\,d\tilde{P}.
\label{eqLemTraNew.pound3}%
\end{equation}
This identification $P\leftrightarrow\tilde{P}$ is based on an implicit
application of Riesz's theorem; see \cite[Chapter 1]{Rud87}.
\end{remark}

\begin{proof}
[Proof of Lemma \textup{\ref{LemTraNew.pound}}]The proof of Kolmogorov's
extension result may be given several forms, but we note that the argument we
used above (in a special case), based on an application of the
Stone--Weierstra\ss \ theorem, also works in general.
\end{proof}

\subsection{\label{ProSpa}The probability space $\Omega$}

We note that the probability space $\Omega$ itself carries mappings $\sigma$
and $\tau_{i}$, $i=0,1,\dots,N-1$.
Stressing the $\Omega$-dependence, we write%
\begin{equation}
\sigma^{\Omega}\left(  \omega\right)  =\left(  \omega_{2},\omega_{3}%
,\dots\right)  \text{\quad and\quad}\tau_{i}^{\Omega}\left(  \omega\right)
=\left(  i,\omega_{1},\omega_{2},\dots\right)  \label{eqTra.15}%
\end{equation}
for $\omega=\left(  \omega_{1},\omega_{2},\dots\right)  \in\Omega$.

\begin{remark}
\label{RemTraNew.poundbis}The connection between the cylinder sets in
\textup{(\ref{eqTra.8})} and the iterated function systems \textup{(}%
IFS\/\textup{)} $\left(  X,\sigma,\tau_{0},\dots,\tau_{N-1}\right)  $ may be
spelled out as follows: the cylinder sets in $\Omega$ generate the
sigma-algebra of measurable subsets of $\Omega$, and similarly the subsets
$\tau_{i_{1}}\cdots\tau_{i_{n}}\left(  X\right)  \subset X$ generate a
sigma-algebra of measurable subsets of $X$. When nothing further is specified,
these will be the sigma-algebras which we refer to when discussing the
measurable functions on $\Omega$ and $X.$ In particular, we will denote by
$M\left(  \Omega\right)  $ and $M\left(  X\right)  $ the respective algebras
of all bounded measurable functions on $\Omega$, respectively $X$.

Note that if $X=\left[  \, 0,1\,\right]  $, and if
\[
\tau_{i}x=\frac{ x+i}{N}, \qquad0\leq i\leq N-1,
\]
then we recover the familiar $N$-adic subintervals:%
\begin{equation}
\tau_{i_{1}}\cdots\tau_{i_{n}}\left(  X\right)  =\left[  \, \frac{i_{1}}%
{N}+\dots+\frac{i_{n}}{N^{n}},\,\frac{i_{1}}{N}+\dots+\frac{i_{n}}{N^{n}%
}+\frac{1}{N^{n}}\,\right]  . \label{eqRemTraNew.poundbis}%
\end{equation}

\end{remark}

\begin{lemma}
\label{LemTraNew.poundpound}There is a unique mapping $\rho\colon M\left(
\Omega\right)  \rightarrow M\left(  X\right)  $ which satisfies%
\begin{equation}
\rho\left(  fg\right)  =\rho\left(  f\right)  \rho\left(  g\right)
\label{eqTraNew.star}%
\end{equation}
and%
\begin{equation}
\rho\left(  \chi_{A\left(  i_{1},\dots,i_{n}\right)  }^{{}}\right)
=\chi_{\tau_{i_{1}}\cdots\tau_{i_{n}}\left(  X\right)  }^{{}}.
\label{eqTraNew.poundprime}%
\end{equation}
The mapping $\rho$ is an isomorphism of $M\left(  \Omega\right)  $ onto
$M\left(  X\right)  $.
\end{lemma}

\begin{proof}
Recalling (\ref{eqTra.14}), we note that%
\[
\chi_{A\left(  i_{1},\dots,i_{n}\right)  }^{{}}\left(  \omega\right)
=\delta_{i_{1},\omega_{1}}\cdots\delta_{i_{n},\omega_{n}},\qquad\omega
\in\Omega.
\]
As a result,%
\begin{equation}
\chi_{A\left(  i_{1},\dots,i_{n}\right)  }^{{}}\chi_{A\left(  j_{1}%
,\dots,j_{n}\right)  }^{{}}=\delta_{i_{1},j_{1}}\cdots\delta_{i_{n},j_{n}}%
\chi_{A\left(  i_{1},\dots,i_{n}\right)  }^{{}}. \label{eqTraNew.pound3}%
\end{equation}
We then define $\rho$ first on $\mathfrak{A}_{n}$ by
\[
\rho\left(  \sum_{i_{1},\dots,i_{n}}a_{i_{1},\dots,i_{n}}\chi_{A\left(
i_{1},\dots,i_{n}\right)  }^{{}}\right)  =\sum_{i_{1},\dots,i_{n}}%
a_{i_{1},\dots,i_{n}}\chi_{\tau_{i_{1}}\cdots\tau_{i_{n}}\left(  X\right)
}^{{}},
\]
where $a_{i_{1},\dots,i_{n}}\in\mathbb{C}$, and note that (\ref{eqTraNew.star}%
) is satisfied.

It is easy to check that the extension of $\rho$ from $\mathfrak{A}_{n}$ to
$\mathfrak{A}_{n+1}$ is consistent. The final extension from $\bigcup
_{n}\mathfrak{A}_{n}$ to $M\left(  \Omega\right)  $ is done by Kolmogorov's
lemma, and it can be checked that $\rho$ has the properties stated in the
conclusion of the present lemma.
\end{proof}

\begin{theorem}
\label{LemTra.3}Let $X$, $W$, and $N$ be as described in the beginning of this
chapter, and let $\left\{  \,P_{x}\mid x\in X\,\right\}  $ be the process
obtained in the conclusion of Lemma \textup{\ref{LemTra.2}}. Then%
\begin{equation}
\sum_{i=0}^{N-1}W\left(  \tau_{i}x\right)  P_{\tau_{i}x}\left[  \, f\left(
i,\,\cdot\,\right)  \,\right]  =P_{x}\left[  \, f\,\right]  \text{\qquad for
all }f\in C\left(  \Omega\right)  . \label{eqTra.16}%
\end{equation}
Moreover, equation \textup{(\ref{eqTra.16})} determines $\left( P_{x}\right)  $ uniquely.
\end{theorem}

\begin{remark}
\label{RemTra.4half}Stated informally, formula \textup{(\ref{eqTra.16})} is an
assertion about the random walk: it says that if the walk starts at $x$, then
with probability one, it makes a transition to one of the $N$ points $\tau
_{0}\left(  x\right)  $, $\dots$, $\tau_{N-1}\left(  x\right)  $. The
probability of the move $x\mapsto\tau_{i}x$ is $W\left(  \tau_{i}x\right)  $.
Recall \textup{(\ref{eqTra.9})} asserts that $\sum_{i}W\left(  \tau
_{i}x\right)  =1$.
\end{remark}

\begin{proof}
[Proof of Theorem \textup{\ref{LemTra.3}}]It follows from (\ref{eqTra.12}) and
the arguments in the proof of Lemma \ref{LemTra.2} that it is enough to verify
(\ref{eqTra.16}) for $f\in C_{\operatorname*{fin}}\left(  \Omega\right)  $, or
for $f\in\mathfrak{A}_{n}$. Let $f\in\mathfrak{A}_{n}$. Then
\begin{align*}
&  \sum_{i}W\left(  \tau_{i}x\right)  P_{\tau_{i}x}\left[  \, f\left(
i,\,\cdot\,\right)  \,\right] \\
&  \qquad=\sum_{i}\;\;\;\;\sum
_{\makebox[0pt]{\hss$\scriptstyle  \omega_{1},\dots,\omega_{n}$\hss}}%
\;W\left(  \tau_{i}x\right)  W\left(  \tau_{\omega_{1}}\tau_{i}x\right)
\cdots W\left(  \tau_{\omega_{n}}\cdots\tau_{\omega_{1}}\tau_{i}x\right)
f\left(  i,\omega_{1},\dots,\omega_{n}\right) \\
&  \qquad\underset
{\makebox[0pt]{\hss$\scriptstyle  \text{by (\ref{eqTra.12})}$\hss}}%
{=}\;\;\;\;P_{x}\left[  \, f\,\right]  ,
\end{align*}
which is the desired conclusion. Recall $C_{\operatorname*{fin}}\left(
\Omega\right)  $ is norm-dense in $C\left(  \Omega\right)  $.

Note that the formula \textup{(\ref{eqTra.16})}
generalizes the familiar notion of selfsimilarity for measures introduced by
Hutchinson in \cite{Hut81}. In fact,
\textup{(\ref{eqTra.16})} may be restated as%
\begin{equation}
\sum_{i=0}^{N-1}W\left(  \tau_{i}x\right)  P_{\tau_{i}x}\circ\left(  \tau
_{i}^{\Omega}\right)  ^{-1}=P_{x}. \label{eqRemTraNew.pound}%
\end{equation}
The uniqueness part of the theorem follows then from Lemma \ref{LemTraNew.pound}.
\end{proof}

\begin{acknowledgements}
The author is pleased to acknowledge numerous constructive
discussions about infinite products with Professors Dorin Dutkay and Richard
Gundy; with Akram Aldroubi,
John Benedetto and other members of our NSF-Focused Research
Group (FRG). And we thank Brian Treadway for outstanding typesetting and
helpful suggestions.
\end{acknowledgements}

\ifx\undefined\bysame
\newcommand{\bysame}{\leavevmode\hbox to3em{\hrulefill}\,}
\fi

\end{document}